\documentclass[graybox]{svmult}

\usepackage{mathptmx}       

\usepackage{helvet}         
\usepackage{courier}        

\usepackage{graphicx}              

\usepackage{siunitx}
\usepackage{amsbsy}
\usepackage{amsmath,bm}
\usepackage{amssymb}
\usepackage{amsfonts}
\usepackage{algorithmic}
\usepackage{algorithm} 

\usepackage{multirow}

\usepackage[bottom]{footmisc}
\usepackage{cite}

\usepackage{url}

\newcommand{\R}{\mathbb R}

\renewcommand{\rho}{\varrho}

\newcommand{\yd}[0]{y_\text{data}}
\newcommand{\J}[0]{\mathcal J}

\begin{document}

\title*{Structure-preserving identification of port-Hamiltonian systems --- a sensitivity-based approach}
\titlerunning{Structure-preserving identification of port-Hamiltonian systems}
\author{Michael Günther 
\and
Birgit Jacob  
\and
Claudia Totzeck 
}
\institute{Michael Günther, Birgit Jacob and Claudia Totzeck \at Bergische Universität Wuppertal, IMACM, Gaußstraße 20, D-42119 Wuppertal,
\email{[guenther,bjacob,totzeck]@uni-wuppertal.de}}

%
%
\maketitle

\abstract*{We present a gradient-based calibration algorithm to identify a port-Hamiltonian system from given input-output data. The gradient is computed with the help of sensitivities and the algorithm is tailored such that the structure of the system matrices of the port-Hamiltonian system (skew-symmetry and positive semi-definitness) is preserved in each iteration of the algorithm. As we only require input-output data, we need to calibrate the initial condition of the internal state of the port-Hamiltonian system as well. Numerical results with synthetic data show the feasibility of the approach.}

\abstract{We present a gradient-based calibration algorithm to identify a port-Hamiltonian system from given time-domain input-output data. The gradient is computed with the help of sensitivities and the algorithm is tailored such that the structure of the system matrices of the port-Hamiltonian system (skew-symmetry and positive semi-definitness) is preserved in each iteration of the algorithm. As we only require input-output data, we need to calibrate the initial condition of the internal state of the port-Hamiltonian system as well. Numerical results with synthetic data show the feasibility of the approach.}
  
  \section{Introduction}
 In structure-preserving modelling of coupled dynamical systems the port-Hamil\-tonian framework allows for constructing overall port-Hamiltonian systems (PHS) provided that (a) all subsystems are PHS and (b) a linear coupling between the input and outputs of the subsystems is provided~\cite{MeMo19, vanDerSchaft06,EbMS07,DuinMacc09}.
 In realistic applications this approach reaches its limits: for a specific subsystem, either no physics-based knowledge is available which allows for defining a physics-based PHS or (b) one is forced to use user-specified simulation packages with no information of the intrinsic dynamics, and thus only the input-output characteristics are available. 
 
 In both cases a remedy for such a subsystem is as follows: generate input-output data either by physical measurements or evaluation of the simulation package, and based on that derive a PHS surrogate that fits these input-output data best. This PHS surrogate can than be used to model the subsystem, and overall one gets a coupled PHS with structure-preserving properties.
 
 Our approach aims at constructing a best-fit PHS model in one step, without the need of first deriving a best-fit linear state-space model and then, in a post-processing step, finding the nearest port-Hamiltonian realization, see, for example, \cite{ChMH19,ChGB22}. In contrast to approaches such as \cite{BGD20} we follow a time domain approach. \cite{SCH} uses a time domain approach as well by  parametrization of the class of PHS is used which permits the usage of unconstrained optimization solvers during identification. Here a PHS with $n$ states and $k$ inputs and outputs is represented by $n(\frac{3n+1}{2}+2k)+k^2$ parameters. In this article, we   develop a gradient-based calibration algorithm to identify a PHS from given time-domain input-output data. The gradient is
computed with the help of sensitivities. 
 
 Consequently, we thus consider the surrogate PHS  system given by
\begin{subequations}\label{totzeck:eq:state}
\begin{align}
    \frac{d}{dt} x &= (J-R) Qx + Bu, \qquad x(0)=\hat x,\\
    y &= B^\top Q x,
\end{align}
\end{subequations}
where $J,Q,R\in \R^{n\times n}$ with $J=-J^\top, Q>0, R \ge 0.$ We assume to have some given reference data $\yd$ and  $B \in \R^{n \times k}$ as well as the input signal $u.$ 

The task is to fit the system matrices and the initial conditions $v=(J,Q,R,\hat x)$ to the data. We therefore define the cost functional
\[
\J(x,v) = \frac{1}{2} \int_0^T |y(t) - \yd(t) |^2 dt = \frac{1}{2} \int_0^T |B^T Q x(t) - \yd(t) |^2 dt
\]
leading us to the calibration problem 
\begin{equation}\label{totzeck:eq:calibration} \tag{P}
    \min \J(x,v) \quad \text{ subject to} \quad \eqref{totzeck:eq:state}.
\end{equation}


As we are only interested in the input-output behaviour of the system, we can eliminate $Q$ from the dynamics. In fact, by Cholesky decomposition we obtain $V$ with $Q = V V^\top$. 
\begin{equation*}
    w=V^\top x,\; \tilde B= V^\top B,\, \tilde J= V^\top JV,\; \tilde R=V^\top  R V
\end{equation*}
yields the system
\begin{align}
    \frac{d}{dt} w &= (\tilde J- \tilde R) w + \tilde Bu, \qquad w(0)=\hat w(=V^\top  \hat x), \label{totzeck:eq:ODE}\\
    y &= \tilde B^\top w. \label{totzeck:eq:output}
\end{align}
For later use we define the state operator $e$ corresponding to \eqref{totzeck:eq:ODE} as
\[
e(w,v) = \begin{pmatrix} \frac{d}{dt} w - (\tilde J- \tilde R) w - \tilde Bu \\ w(0) - w_0\end{pmatrix}.
\]
Hence, \eqref{totzeck:eq:ODE} is equivalent to $e(w,v)=0.$

The transformed cost functional is given by
\[
\tilde J(w,v) = \frac{1}{2} \int_0^T |y(t) - \yd(t) |^2 dt = \frac{1}{2} \int_0^T |\tilde B^T w(t) - \yd(t) |^2 dt.
\]

After the transformation we are left to identify the matrices $\tilde J, \tilde R$ and $w_0.$ For notational convenience we define the space of admissible controls
$$\mathcal V = \{ (\tilde J, \tilde R,w_0) \in \R^{n\times n} \times \R^{n\times n} \times \R^n\, \colon \, \tilde J^\top = -\tilde J,\; \tilde R \ge 0 \}.$$

Note that the system of differential equations admits a unique solution by standard ODE theory. This allows us to define the control to state map $$S \colon \mathcal V \mapsto C([0,T],\R^n), \quad S(v) = w. $$
Moreover, we use $S$ to define the reduced cost functional
\[
\hat J(v) := \frac{1}{2} \int_0^T |\tilde B^T S(v)(t) - \yd(t) |^2 dt.
\]
In the following we aim to derive an gradient-based algorithm that allows us to solve the calibration problem numerically. In particular, we require to compute the gradient of $\hat J.$ Details are presented in the next section. From now on we only work with the transformed system and drop the $\sim$ for notational convenience.

\section{Sensitivity approach}
We emphasize that the system matrices $J,R$ as well as the initial condition $\hat x$ are finite dimensional. It is therefore feasible to employ an sensitivity approach \cite{HPUU} for the calibration problem. 

To compute the sensitivities require admissible directions for the G\^ateaux derivatives. Due to the structural restrictions, $J$ can only be varied in direction $h_J$ satisfying $h_J^\top = -h_J$ and $R$ can only be varied by symmetric matrices. 

The directional derivative of $\hat J$ in direction $h = (h_J, h_R, h_x)$ is given by
\[
d\hat J(v)[h] = \langle \hat J'(v), h \rangle = \langle d_w J(w,v), S'(v)h \rangle + \langle d_v J(w,v), h \rangle
\]
 To evaluate this, we require $dw(v,h)= S'(v)h$ the so-called sensitivity. Here, we make use of the state equation $e(w,v) = 0.$ In fact, it holds
\begin{equation}\label{totzeck:eq:sensitivity}
e_w(w,v) dw(v,h) + e_v(w,v) = 0 \qquad \Leftrightarrow \qquad  e_w(w,v) dw(v,h) = -e_v(w,v)h.
\end{equation}
We emphasize that in order to identify the gradient $\hat J'(v)$ we need to compute the directional derivative w.r.t.~all basis element of the tangent space of $\mathcal V.$

\section{Gradient-descent algorithm}
In the previous section we established the theoretical foundation of the gradient descent algorithm we present in the following.

Starting from an initial guess of system matrices and initial condition $v_0 = (J_0, R_0, \hat x_0)$ we compute the sensitivities $dw(v,h)$ for all basis elements of the tangent space of $\mathcal V$ by solving \eqref{totzeck:eq:sensitivity} and use the sensitivity information to evaluate the gradient $\hat J'(v_0).$ Then we seek for an admissible stepsize $\sigma$ using Armijo-rule \cite{HPUU}, see the pseudo code in Algorithm~\ref{totzeck:alg:armijo} and update the system matrices and the initial condition $v_0 \gets v_0 - \sigma \hat J'(v_0).$  The calibration procedure is stopped when the cost functional value is sufficiently small. A pseudo code of the calibration algorithm can be found in Algorithm~\ref{totzeck:alg:calibration}.\\

\begin{algorithm}[ht!]
\caption{Armijo step size search}
\begin{algorithmic}\label{totzeck:alg:armijo} 
\REQUIRE gradient $g,$ initial step size $\sigma$ and safety parameter $\gamma$ 
\ENSURE admissible step size $\sigma,$ new parameter set $v'$
\STATE $v' \gets v + \sigma g$
\WHILE{$\hat J(v') - \hat J(v') > -\gamma \sigma \| g \|^2 $}
\STATE $\sigma \gets 0.5 \sigma$
\STATE $v' \gets v- \sigma g$
\ENDWHILE
\end{algorithmic}
\end{algorithm}

\begin{algorithm}[ht!]
\caption{Gradient-based calibration algorithm}
\begin{algorithmic}\label{totzeck:alg:calibration} 
\REQUIRE initial guess $v_0$ and additional parameters
\ENSURE calibrated system matrices and initial condition $v=(J,R,\hat x)$
\WHILE{$\hat J(v_0) > \epsilon_\text{stop}$}
\FOR{all admissible directions $h$}
\STATE compute $dw(v_0,h)$ by solving \eqref{totzeck:eq:sensitivity}
\ENDFOR
\STATE identify $\hat J'(v_0)$ 
\STATE find admissible step size $\sigma$ by Armijo-rule, see Algorithm~\ref{totzeck:alg:armijo}
\STATE $v_0 \gets v_0 - \sigma \hat J'(v_0)$
\ENDWHILE
\end{algorithmic}
\end{algorithm}

The presented algorithm can be used for numerical studies. In the following we discuss a proof of concept with states $x \in C([0,T],\R^2).$ 

\section{Proof of concept}
In the following we discuss a proof of concept with states $x \in C([0,T],\R^2)$ and output $y \in C([0,T],\R).$ In the two dimensinal setting the basis elements of the tangent space of $\mathcal V$ are manageble. Indeed, we have the basis elements
\[
J_1 = \begin{pmatrix} 0 & -1 \\ 1 & 0 \end{pmatrix},\, R_1  = \begin{pmatrix} 1 & 0 \\ 0 & 0 \end{pmatrix},\, R_2  = \begin{pmatrix} 0 & 0 \\ 0 & 1 \end{pmatrix},\, R_3  = \begin{pmatrix} 0 & 1 \\ 1 & 0 \end{pmatrix},\, x_1 =  \begin{pmatrix} 1 \\ 0 \end{pmatrix},\, x_2 =\begin{pmatrix} 0 \\ 1 \end{pmatrix}.
\]
We assume that $B = \begin{pmatrix} 1 & 1 \end{pmatrix}$ is known and that input signals at the time steps $t_k$ are given as $u(t_k) = 1 + 0.1 N(0,1)$ where $N(0,1)$ denotes a realization of a normally distributed random variable with mean $0$ and standard deviation $1.$ 

For simplicity, we assume that the time steps $t_k, k=1,\dots,K$ coincide with the time step of the Euler discretization that is implemented to solve the state ODE. Indeed, with the initial guess we solve \eqref{totzeck:eq:ODE} using the Euler scheme. Then we obtain the output $y$ by \eqref{totzeck:eq:output}, which we use to evaluate the cost functional for the initial guess. If the cost values is higher than the tolerance $\epsilon_\text{stop}$ we start the calibration procedure. 

For notational convenience we split the sensitivity $dw(v,h)$ into the parts $h_J, h_R$ and $h_x.$ The sensitivity w.r.t.~$J$ is computed by solving $e_w(w,v) dw(v,h_J) = -e_v(w,v)h_J$ which can be written explicitly as
\[
\frac{d}{dt} dw(v,h_J) - (J-R) dw(v,h_J) = h_Jw, \qquad dw(v,h_J)(0) = 0. 
\]
In the two dimensional case, there is only one admissible direction $h_J = J_1.$ For the sensitivities w.r.t.~$R$ we solve
\[
\frac{d}{dt} dw(v,h_R) - (J-R) dw(v,h_R) = -h_Rw, \qquad dw(v,h_R)(0) = 0
\]
for $h_R = \{R_1,R_2,R_3\}.$ For the initial condition we solve
\[
\frac{d}{dt} dw(v,h_x) - (J-R) dw(v,h_x) = 0, \qquad dw(v,h_x)(0) = h_x 
\]
for $h_x = \{x_1, x_2\}.$ 

The directional derivative of the cost functional reads
\[
d\hat J(v)[h] = \langle B^\top S(v) - y_\text{data}, B^\top S'(v) h \rangle = \int_0^T B\big(B^\top S(v)(t) - y_\text{data}(t)\big),  (S'(v) h)(t) dt,
\]
which we can evaluate with the help of the sensitivities computed above. Note that $d\hat J(v)[h_\bullet] \in \R$ for all $h_\bullet$ discussed above. Hence, the gradient is assembled as follows
\[
\hat J'(v) = \left[ d\hat J(v)[J_1] J_1 \quad \sum_{\ell=1}^3 d\hat J(v)[R_\ell] R_\ell \quad\sum_{\ell=1}^2 d\hat J(v)[x_\ell] x_\ell \right]^\top
\]

\section{Numerical results}
\begin{figure}[ht!]
\sidecaption
\includegraphics[scale=.4]{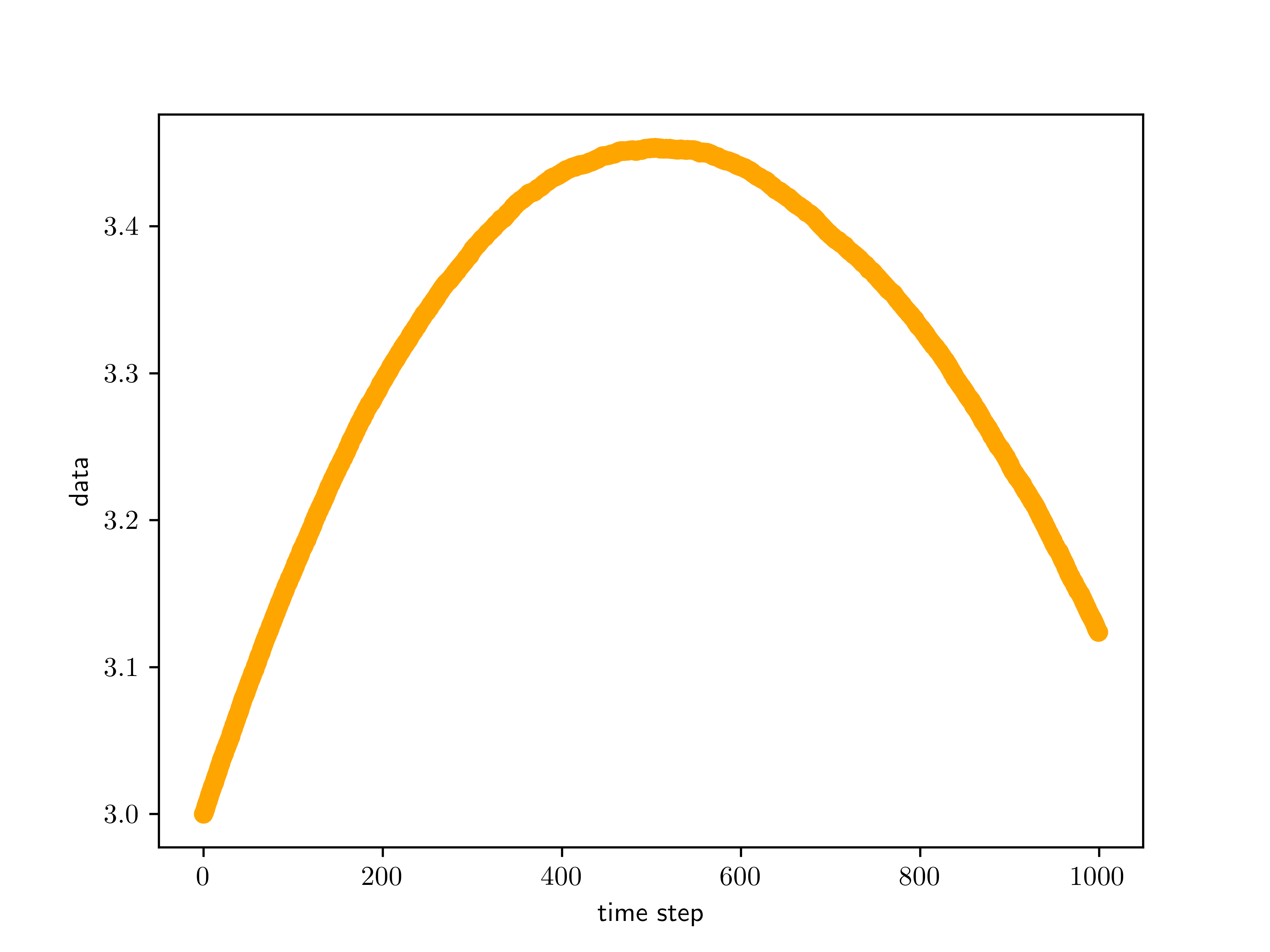}
\includegraphics[scale=.4]{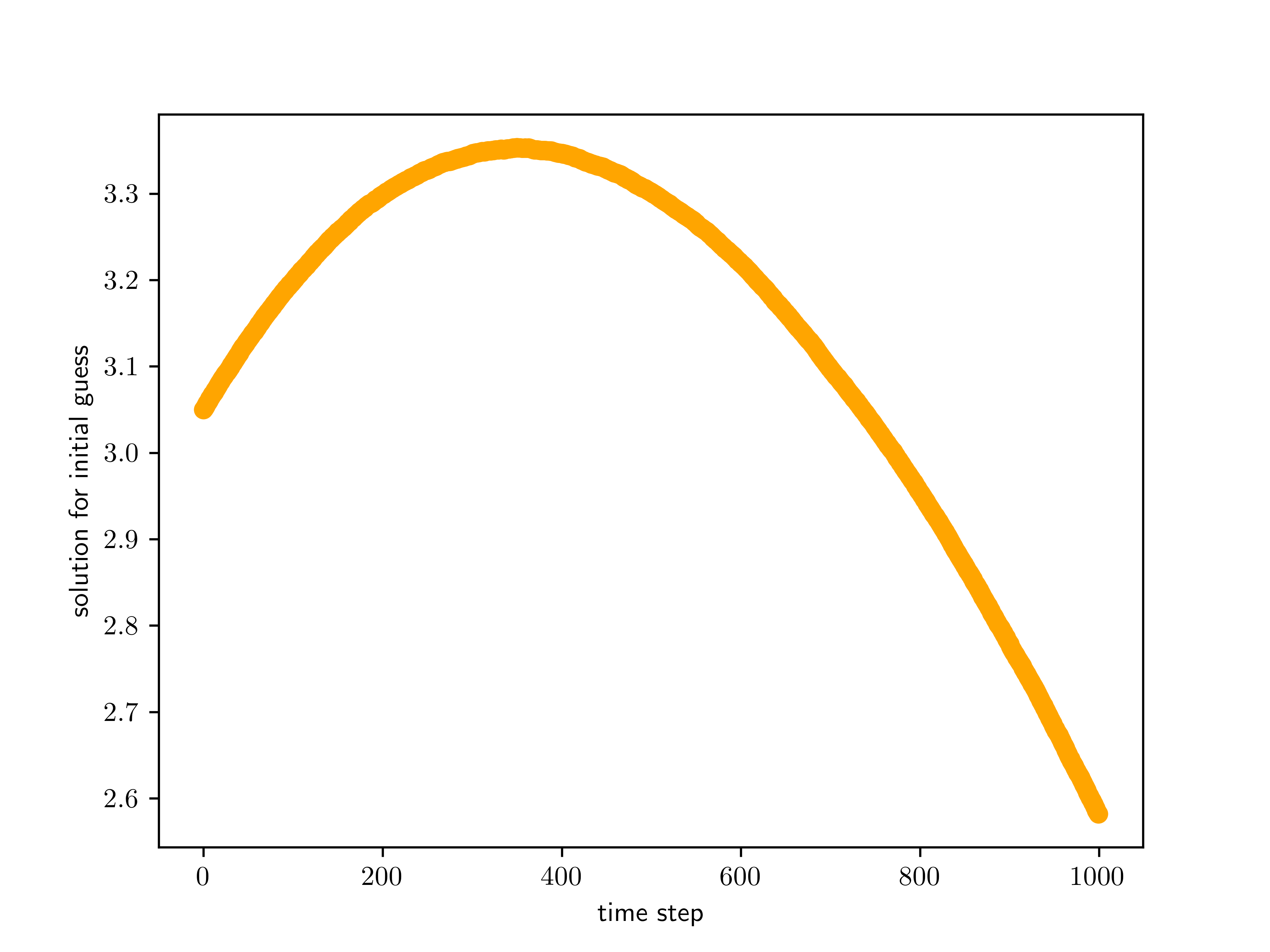}
\caption{Left: output $y_\text{data}$ corresponding to the data given in \eqref{totzeck:eq:data}. Right: output $y_0$ corresponding to the initial guess \eqref{totzeck:eq:guess}.}
\label{totzeck:fig:ydata}       
\end{figure}

For our proof of concept we generate synthetic data by solving the state system for fixed data matrices $J_\text{data}, R_\text{data}$ and initial condition $\hat x_\text{data}.$ For the following results we choose
\begin{equation}\label{totzeck:eq:data}
J_\text{data} = \begin{pmatrix} 0 & 1 \\ -1 & 0\end{pmatrix},\quad R_\text{data} = \begin{pmatrix} 0.5 & 0 \\ 0 & 0.3  \end{pmatrix}, \quad\hat x_\text{data} = \begin{pmatrix} 1 \\ 2  \end{pmatrix}.
\end{equation}
The data yields the reference output $y_\text{data}$ shown in Figure~\ref{totzeck:fig:ydata} (left).

We start the proof of concept with the initial guess given by
\begin{equation}\label{totzeck:eq:guess}
J_0 = \begin{pmatrix} 0 & 1.2 \\ -1.2 & 0\end{pmatrix},\quad R_0 = \begin{pmatrix} 0.4 & 0 \\ 0 & 0.4  \end{pmatrix}, \quad\hat x_0 = \begin{pmatrix} 1.1 \\ 1.95  \end{pmatrix}
\end{equation}
leading to the output in Figure~\ref{totzeck:fig:ydata} (right). We set $T=1$ and use $1000$ time steps for the Euler discretization. The Armijo-search for an admissible step size is initialized with $\sigma = 10$ and the $\sigma \gets \sigma/2$ if the current step size is not admissible.

\begin{figure}[ht!]
\sidecaption
\includegraphics[scale=.4]{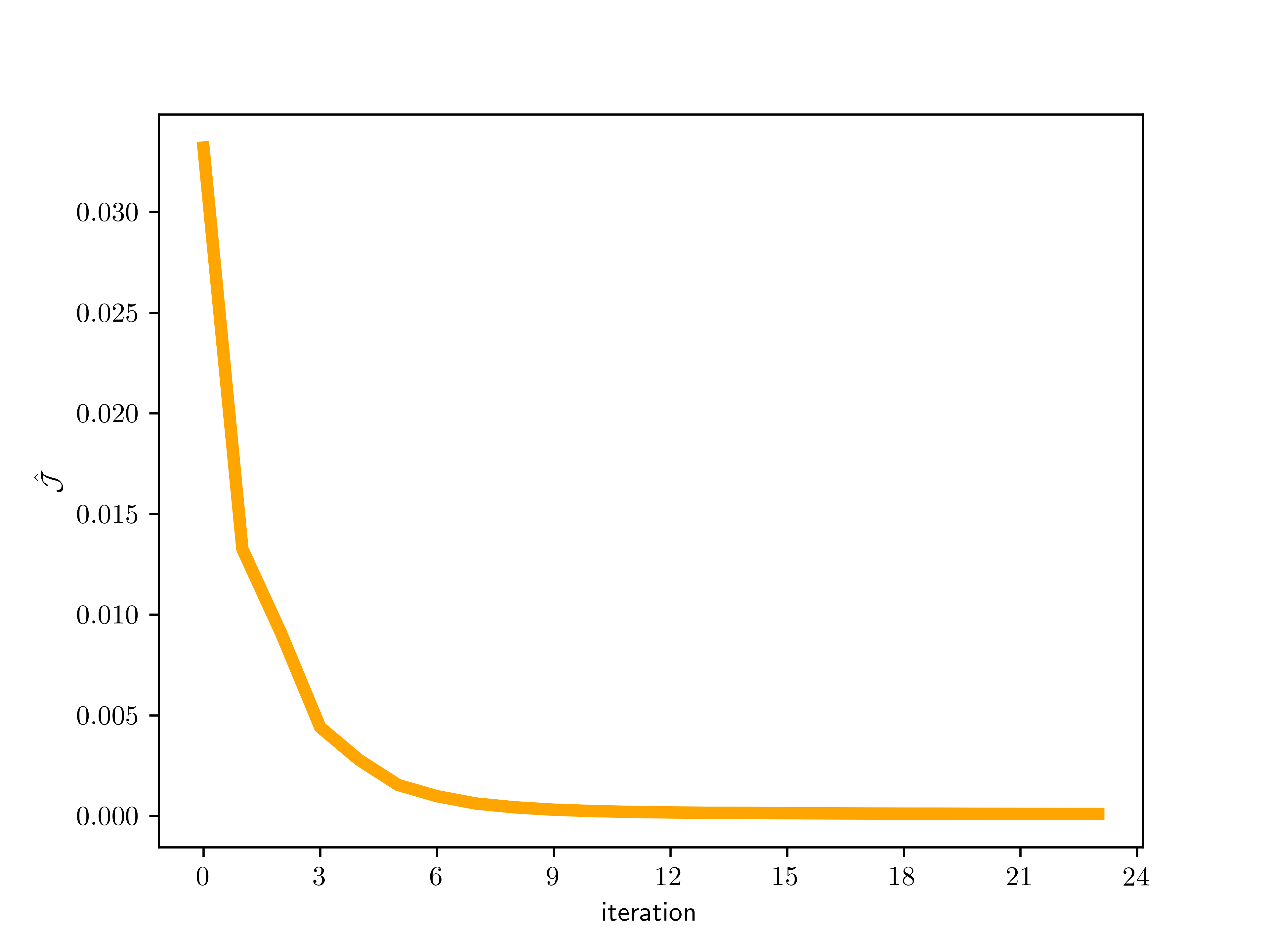}
\includegraphics[scale=.4]{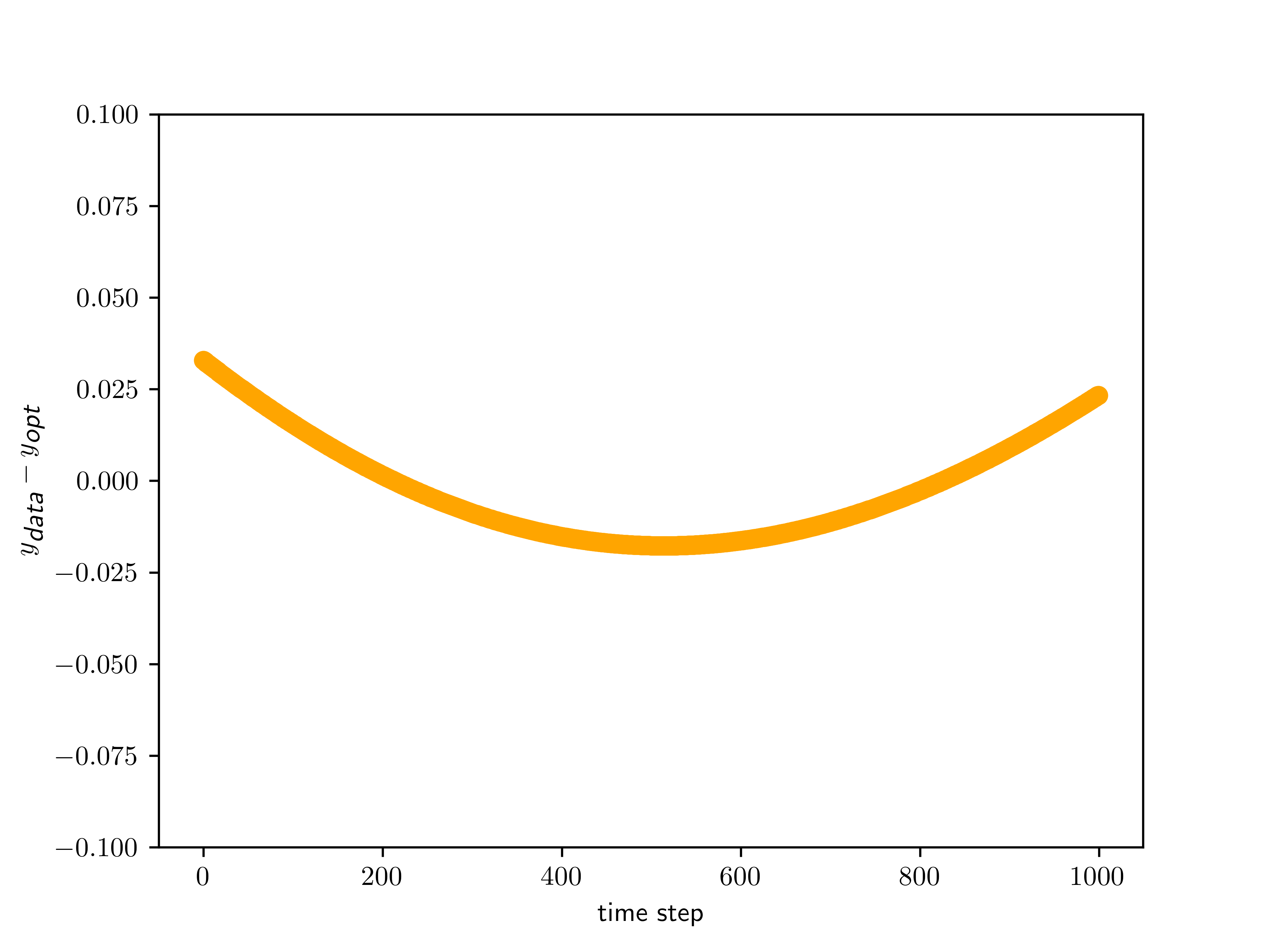}
\caption{Left: output $y_\text{data}$ corresponding to the data given in \eqref{totzeck:eq:data}. Right: difference of reference output and output of the calibrated system.}
\label{totzeck:fig:cost}       
\end{figure}

\begin{figure}[ht!]
\sidecaption
\includegraphics[scale=.4]{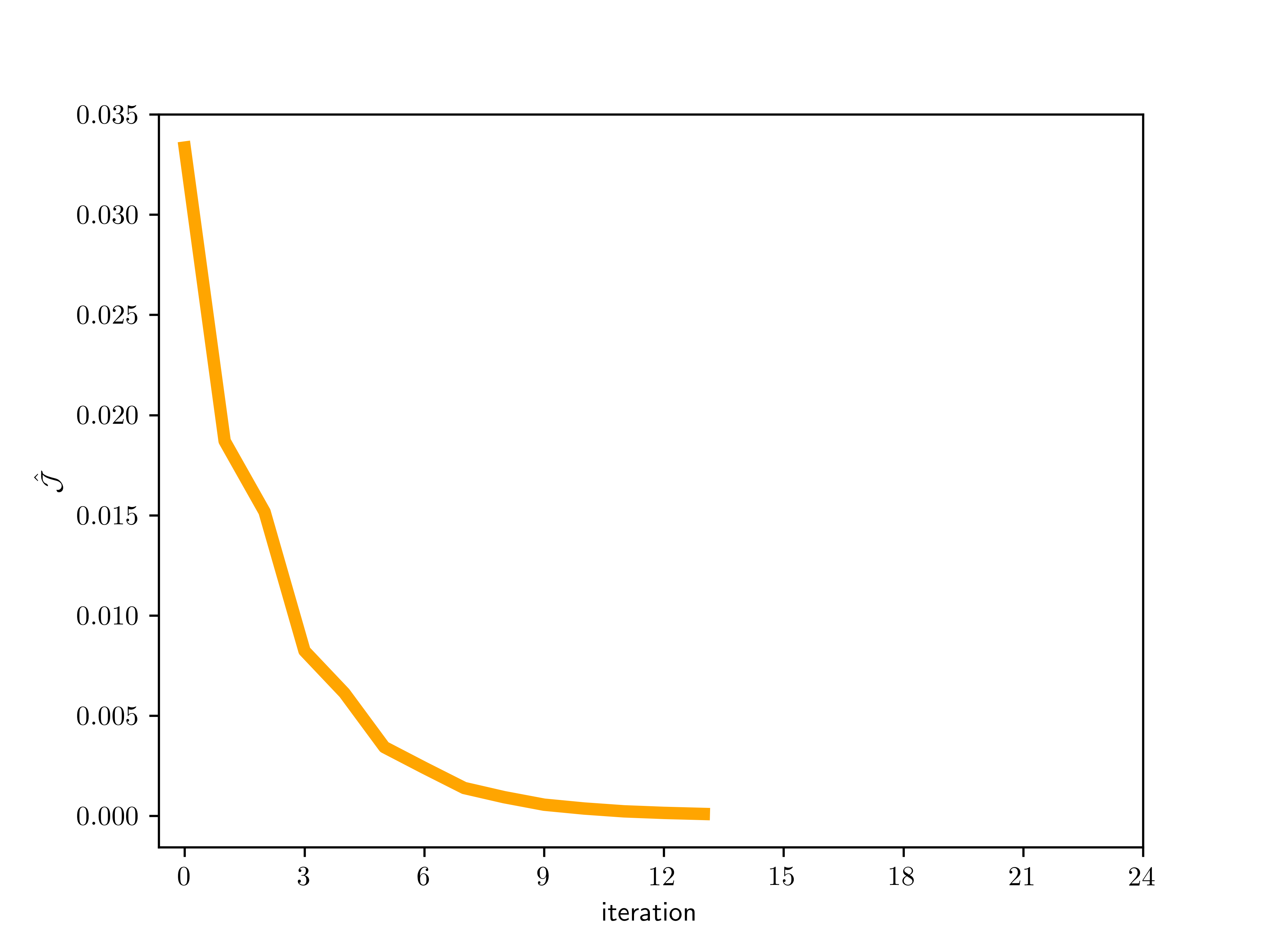}
\includegraphics[scale=.4]{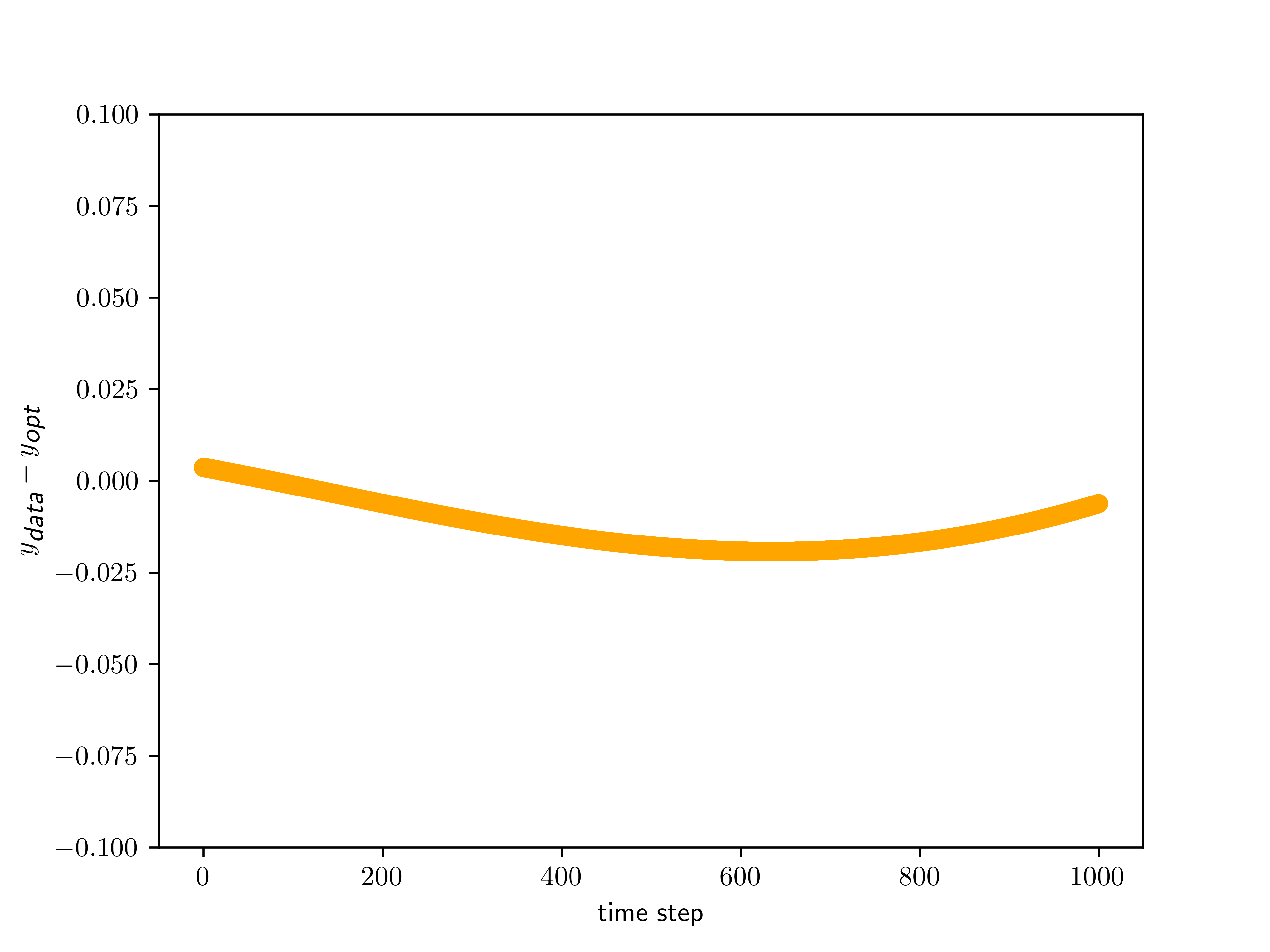}
\caption{Left: output $y_\text{data}$ corresponding to the data given in \eqref{totzeck:eq:data}. Right: difference of reference output and output of the calibrated system with $R$ diagonal.}
\label{totzeck:fig:cost2}       
\end{figure}

Algorithm~\ref{totzeck:alg:calibration} is able to reproduce the output $y_\text{data}$ with $\epsilon_\text{stop} = 1e^{-4}$ in $22$ gradient steps. The evolution of the cost function is shown in Figure~\ref{totzeck:fig:cost} (left) and the difference $y_\text{data} - y_\text{opt}$ is plotted in Figure~\ref{totzeck:fig:cost} (right).
The calibrated matrices and initial data read
\[
J_\text{opt} = \begin{pmatrix} 0 & 1.073 \\ -1.073 & 0\end{pmatrix},\quad R_\text{opt} = \begin{pmatrix} 0.379 & -0.080 \\ -0.080 & 0.367  \end{pmatrix}, \quad\hat x_\text{opt} = \begin{pmatrix} 1.039 \\ 1.929  \end{pmatrix},
\]
where we rounded to precision $1e^{-3}.$ It jumps to the eye that $R_\text{opt}$ has nonzero off-diagonal entries. Out of curiosity we run the same toy problem with $R$ restricted diagonal matrices. We obtain the calibrated matrices and initial data 
by
\begin{equation}
    J_\text{opt,2} = \begin{pmatrix} 0 & 1.016 \\ -1.016 & 0\end{pmatrix},\quad R_\text{opt,2} = \begin{pmatrix} 0.351 & 0 \\ 0 & 0.335  \end{pmatrix}, \quad\hat x_\text{opt,2} = \begin{pmatrix} 1.023 \\ 1.973  \end{pmatrix},
\end{equation}
again rounded to precision $1e^{-3}.$ The additional structural information in $R$ yields overall to better calibrated results. Compare Figure~\ref{totzeck:fig:cost2} for the cost evolution and the difference of the outputs for the calibration with $R$ restricted to diagonal matrices.

\section{Conclusion and outlook}
We present a gradient-based algorithm to identify a port-Hamiltonian system consisting of ordinary differential equation to given input-output data. The gradient is computed with the help of a sensitivity approach. A proof of concept shows the feasibility of the approach. 

As the effort of the sensitivity approach scales with the number of basis elements of the tangent space, the proposed calibration algorithm is only recommended for small systems. In future work, we investigate an adjoint-based approach to compute the gradient in order to derive a structure-preserving calibration algorithm for port-Hamiltonian input-output systems.



\end{document}